\documentclass[centertags,leqno]{article}

\usepackage{amsmath}
\usepackage{amssymb}
\usepackage{latexsym}
\usepackage{amsxtra}
\usepackage{amscd}
\usepackage{theorem}
\usepackage{graphics}

\usepackage{epic}

\theoremstyle{plain}

{\theorembodyfont{\slshape}

        \newtheorem{thm}{Theorem}[section]
        
        \newtheorem{lemma}[thm]{Lemma}

}

{\theorembodyfont{\rmfamily}

        \newtheorem{defn}[thm]{Definition}
        
        \newtheorem{rem}[thm]{Remark}

}

\renewcommand{\em}{\sl}

\newcommand{\proof}{\noindent {\bf Proof:\ }}
\newcommand{\Endproof}{\hspace*{\fill} $\Box$ \vspace{1ex} \noindent }

\def\rk{{\rm rk}}
\def\F{{\cal F}}

\def\M{{\bf M}}

\newcommand{\Mod}{ {\rm Mod} }

\def\K{{\cal K}}
\newcommand{\CC}{{\Bbb C}}

 \def\PP{{\Bbb P}}

\def\GL{{\rm GL}}
\def\SL{{\rm SL}}

\def\PVI{{\rm PVI}}

\def\L{{\cal L}}

\def\A{{\bf A}}

\def\J{{\bf J}}
\def\K{{\cal K}}

\newcommand{\CCC}{{\cal  C}}
\newcommand{\NNN}{{\cal  N}}
\newcommand{\Hom}{{\rm  Hom}}
\newcommand{\Tr}{{\rm  Tr}}

\title{Painlev\'e equations and the middle convolution}

\author{
   Michael Dettweiler$^*$ \and
   Stefan Reiter\thanks{both authors gratefully acknowledge
financial support from the  Research Training Network (Galois
Theory and Explicit Methods in Arithmetic) of  the European
Community.}\\{}}

\begin{document}

\maketitle

\begin{abstract} We use the middle convolution to obtain some old and 
new algebraic solutions of the Painlev\'e VI equations. 
\end{abstract}

\section{Introduction}

A description of all  irreducible and physically rigid local
systems on the punctured affine line was given by Katz
\cite{Katz97}. The main tool here-fore is a middle convolution
functor on the category of  perverse sheaves, loc. cit., Chap.~5.
This functor is denoted by $MC_\chi,$ for $\chi$ a one-dimensional
representation of $\pi_1({\Bbb G}_m).$ It
  preserves
important properties of local systems (resp. perverse sheaves) 
like the index of rigidity
and irreducibility,  but in general, $MC_\chi$ changes the rank
and the monodromy group. \\

In \cite{DR}, the authors give a  purely algebraic analogon
 of the functor $MC_\chi$
 (the construction is reviewed in Section \ref{mclambda}).
This  functor is a functor of the category
 of modules of the free group $F_r$ on $r$
generators to itself.  It depends on a scalar $\lambda\in
\CC^\times$ and is denoted by $MC_\lambda.$ Clearly, 
 $MC_\lambda$ can be viewed
as a transformation which sends an $r$-tuple
$(M_r,\ldots,M_1)\in \GL_m^r$ to another 
$r$-tuple $MC_\lambda(M_r,\ldots,M_1)\in \GL_n^r,$
where usually $m\not= n.$ It is shown in  \cite{DR},
that this transformation commutes with the action of the Artin braid
group (up to overall conjugation in $\GL_n$). It is this 
property, which makes the middle convolution functor $MC_\lambda$ 
useful for the study of the Painlev\'e equations.\\

The sixth Painlev\'e equation $\PVI_{\alpha,\beta,\gamma,\delta}$
is a nonlinear  differential 
equation, depending on $4$ parameters
$\alpha,\beta,\gamma,\delta\in \CC:$ 
\begin{eqnarray} y''&=&\frac{1}{2}\left(
\frac{1}{y}+\frac{1}{y-1}+\frac{1}{y-t}\right)(y')^2
-\left(\frac{1}{t}+\frac{1}{t-1}+\frac{1}{y-t}\right)y'\nonumber \\
&&+\; \frac{y(y-1)(y-t)}{t^2(t-1)^2}\left(\alpha+\frac{\beta t}{y^2}
+\frac{\gamma(t-1)}{(y-1)^2}+\frac{\delta t(t-1)}{(y-t)^2}\right).\nonumber
\end{eqnarray}

It is well known, that $\PVI$ arises from isomonodromic deformations
(Schlesinger equations), see \cite{JimboMiwa} and \cite{Hitchin97}.
This approach is recalled in  
 Section \ref{painleve}. Using this 
approach, one can show that finite braid group orbits of triples 
in $\SL_2(\CC)$ give rise to solutions of $\PVI$ with finite branching.
It is expected, that these solutions are always algebraic, see 
\cite{Boalch03}, P. 1012.\\

In this paper, we give a convolution approach 
to the 
 ``Klein solution'' of $\PVI,$ studied by Boalch \cite{Boalch03},
\cite{Boalch04}. Originally, this solution is constructed using 
a braid invariant map  $\varphi$ from triples of pseudo reflections in 
$\GL_3(\CC)$ to triples in $\SL_2(\CC).$ This map is applied to 
a generating triple of reflections of the complex reflection group associated to Klein's simple group of order $168.$ Since the order of this group 
is finite, the generating triple gives rise to a finite braid orbit, 
from which the existence of an algebraic solution of $\PVI$ is derived.

We show that the transformation $\varphi$ is obtained 
from a suitable middle convolution and scaling. Thus a new approach 
to the Klein solution is given. Then it is shown, how to obtain 
new algebraic solutions for $\PVI$ by starting with tuples in 
$\GL_2(\CC).$ Moreover, we give an example of an algebraic solution 
which arises from a rational pullback 
but  cannot be obtained from the middle convolution. One may ask whether 
all algebraic solutions of $\PVI$ arise either by middle convolution and 
scaling or from a rational pullback.\\

The organization of the paper is as follows: In Section 
\ref{mclambda}, we recall the definition and the 
basic facts of $MC_\lambda.$ Section \ref{painleve} is concerned 
with some well known facts about Painlev\'e equations, isomonodromic
deformations and braid orbits. Section \ref{boalch} recalls 
Boalch's map $\varphi$ and gives the connection to $MC_\lambda.$
Finally, in Section \ref{examples} we give some new examples of 
algebraic solutions of $P_{VI},$ using the convolution.

\section{The middle convolution functor $MC_\lambda$}\label{mclambda}

In
this section, we recall the algebraic construction of  the
 convolution functor $MC_\lambda,$ defined in
\cite{DR} and \cite{dr03}. \\

We will use the following notations and conventions throughout the
paper: Let $K$ be
 a field and $G$ a group. The category
of finite dimensional left-$G$-modules is denoted by $\Mod(K[G]).$
Mostly, we do not distinguish notationally between an  element of
$\Mod(K[G])$ and its underlying vector space. Let  $F_r$ denote  the free group on $r$ generators
$f_1,\ldots,f_r.$
 An element
in $\Mod(K[F_r])$ is viewed as a pair $(\M,V),$ where $V$ is a
vector space over $K$ and $\M=(M_r,\ldots,M_1)$ is an element of
$\GL(V)^r$ such that $f_i$ acts on $V$ via $M_i,\, i=1,\ldots,r.$
For $(\M,V)\in \Mod(K[F_r]),$ where $\M=(M_r,\ldots,M_1)\in
\GL(V)^r,$ and $\lambda \in K^\times$ one can construct an element
$(C_\lambda(\M),V^r)\in \Mod(K[F_r]),$
$C_\lambda(\M)=(N_r,\ldots,N_1)\in \GL(V^r)^r,$ as follows: For
$k=1,\ldots,r,$ $N_k$ maps a vector $(v_1,\ldots,v_r)^{\rm tr}$
$\in V^r$ to
\[ \left( \begin{array}{ccccccccc}
                  1 & 0 &  & \ldots& & 0\\
                   & \ddots &  & & &\\
                    & & 1 &&&\\
               (M_1-1) & \ldots&  (M_{k-1}-1)  & \lambda M_{k} &\lambda (M_{k+1}-1) & \ldots
&  \lambda (M_r-1) \\
     &&&&1&&\\
               &   &  & && \ddots  &   \\
                   0 &  &  & \ldots& &0 & 1
          \end{array} \right)\left(\begin{array}{c}
v_1\\
\vdots\\
\\\vdots\\
\\\vdots\\
\\ v_r\end{array}\right)
.\]

 We set $C_\lambda(\M):=(N_r,\ldots,N_1).$
There are the following $\langle N_1,\ldots,N_r \rangle$-invariant
subspaces of  $V^r:$

\[ \K_k = \left( \begin{array}{c}
          0 \\
          \vdots \\
          0 \\
          \ker(M_k-1) \\
            0\\
           \vdots \\
           0
        \end{array} \right)  \quad \mbox{({\it k}-th entry)},\, k=1,\dots,r,\]
and
\[     \L=\cap_{k=1}^r \ker (N_k-1)={\rm ker}(N_1\cdots N_r - 1).
\]
Let $\K:=\oplus_{i=1}^r\K_i.$

\begin{defn}{\rm Let
$MC_\lambda(\M):=(\tilde{N}_r,\dots,\tilde{N}_1)\in
\GL(V^r/(\K+\L))^r,$ where $\tilde{N}_k$ is induced by the action
of $N_k$ on $V^r/(\K+\L).$ The $K[F_r]$-module
$MC_{\lambda}(V):=(MC_\lambda(\A),V^r/(\K+\L))$
 is called the
{\em  middle convolution} of $\M$ with
$\lambda.$}
\end{defn}

\begin{thm} \label{eigen} Let $V=(\M,V)\in \Mod(K[F_r]),$ where
$\M=(M_r,\dots,M_1) \in \GL(V)^r$
and $\lambda \in K^\times.$ Suppose that
$V$ has no $1$-dimensional factors and/or submodules with the
property
 that only one (or none)
of the $M_i$ act non-trivially.\begin{enumerate}

\item If $\lambda \neq 1,$ then $${\rm dim}(MC_\lambda(V))=
\sum_{k=1}^{r} \rk (M_k-1)-
 ({\rm dim}(V)-\rk(\lambda\cdot M_1\ldots M_r-1)).$$

\item If $\lambda_1,\, \lambda_2\in K^\times$ such that
$\lambda_1\lambda_2=\lambda$ and
 $(*)$ and $(**)$ hold for $V,$ then
    \[ MC_{\lambda_2}MC_{\lambda_1}(V)\cong MC_{\lambda}(V).\]

\item Under the assumptions of (ii), if $V$ is irreducible, then
$MC_{\lambda}(V)$ is irreducible.\\

\item   Let   ${\cal B}_r=\langle \beta_1,\ldots,\beta_{r-1} \rangle $ be
the abstract Artin braid group, where  the generators $\beta_1,
\ldots,\beta_{r-1}$ of ${\cal B}_r$ satisfy the usual 
braid relations and act in the following way on tuples
$(g_1,\ldots,g_r)\in G^r$ (where $G$ is a group):
$$
\beta_i(g_r,\ldots,g_1)=(g_r,\ldots,g_{i-1},g_ig_{i+1}g_{i}^{-1},g_i,g_{i+2},
\ldots,g_1),\,\, i=1,\ldots,r-1.$$
 For any $\beta \in{\cal
B}_r$ there exists a $B\in \GL(V^r/(\K+\L))$ such that
$$MC_\lambda(\beta(\M))= \beta(MC_\lambda(\M))^B, $$
where $B$ acts via component-wise conjugation.\\
\end{enumerate}
\end{thm}

\proof (i)-(iv) follow
 analogously to \cite{DR}, Lemma 2.7, Lemma A.4, Theorem 3.5, Corollary 3.6 and
Theorem 5.1 (in this order). \Endproof

A Jordan block of eigenvalue $\alpha\,\in \bar{K}$ and of length
$l$ is denoted by $\J(\alpha,l).$ The following Lemma is a
consequence of \cite{DR}, Lemma 4.1:

\begin{lemma} \label{lemmonodromy1} Let $V=(\M,V)\in \Mod(K[F_r]),$ where
$\M=(M_1,\dots,M_r) \in \GL(V)^r$ such that $\M$ satisfies $(*)$
and $(**).$ Let $\lambda \in K^\times$ and
$MC_{\lambda}(\M)=(\tilde{N}_1,\dots,\tilde{N}_r).$
\begin{enumerate}

\item Every  Jordan
block $\J(\alpha,l)$ occurring in the Jordan
decomposition of $M_i$ contributes a Jordan block $\J(\alpha
\lambda,l')$ to the Jordan decomposition of $\tilde{N}_i,$ where
$$ l':\;=\quad
  \begin{cases}
    \quad l, &
                              \quad\text{\rm if $\alpha \not= 1,\lambda^{-1}$,} \\
    \quad  l-1& \quad \text{\rm if $\alpha =1$,} \\
    \quad l+1, & \quad \text{\rm if $\alpha =\lambda^{-1}$.}
  \end{cases}
  $$
  The only other Jordan blocks which occur in the Jordan
  decomposition of $\
  \tilde{N}_i$ are blocks of the form $\J(1,1).$

\item
Every  Jordan block $\J(\alpha,l)$ occurring in the Jordan
decomposition of  $M_{r+1}$ contributes a Jordan block $\J(\alpha
\lambda,l')$ to the Jordan decomposition of $\tilde{N}_{r+1},$
where
$$ l':\;=\quad
  \begin{cases}
    \quad l, &
                              \quad\text{\rm if $\alpha \not= 1,\lambda^{-1}$,} \\
    \quad  l+1& \quad \text{\rm if $\alpha =1$,} \\
    \quad l-1, & \quad \text{\rm if $\alpha =\lambda^{-1}$.}
  \end{cases}
  $$
  The only other Jordan blocks which occur in the Jordan
  decomposition of $\tilde{N}_{r+1}$ are blocks of the form $\J(\lambda,1).$
\end{enumerate}\end{lemma}

\begin{rem} A geometric interpretation of the middle convolution
can be given as follows, see \cite{dr03}:  
Let $X=\CC\setminus \{a_1,\ldots,a_r\},$
 $$E=\{(x,y)\in \CC^2 \mid x,y\not= a_i,\, i=1,\ldots,r,\, x\not= y\},$$
 ${\rm p}_i:E\to X,\, i=1,2,$ be the $i$-th projection, $$q:E\to \CC^\times,\,
(x,y)\mapsto y-x,$$ $j: E \to \PP^1(\CC)\times X$ the tautological 
inclusion
 and 
$\bar{p}_2: \PP^1(\CC)\times X \to X$ the  (second) projection onto $X.$ 
Moreover, let  $\L_\lambda$ denote the  Kummer sheaf
associated to the representation, which sends 
a generator of $\pi_1(\CC^\times)$ to $\lambda$. 
 Let $\F$ be a 
local system on $X,$ corresponding to 
a tuple $\M:=(M_r,\ldots,M_1)\in \GL_m(\CC)^r$
(by a choice of a homotopy base)
 and let $\lambda \in \CC^\times\setminus 1.$
Then, under the assumptions of Thm. \ref{eigen},  
$ MC_\lambda(\M)$  corresponds to the higher 
direct image sheaf 
$$ R^1 (\bar{p}_2)_*(j_*({\rm p}_1^*(\F)\otimes 
q^*(\L_\lambda)))$$ (under the same choice of homotopy base as before).  
\end{rem}

\section{Painlev\'e VI and the braid group orbits}\label{painleve}

It is well known, that $\PVI$ arises from isomonodromic deformations as 
follows, compare to \cite{JimboMiwa} and \cite{Hitchin97}: Consider 
Fuchsian systems
\begin{equation}\label{fuchs} \frac{d\Phi}{dz}=\left( \sum_{i=1}^3\frac{A_i}{z-a_i}\right)\Phi\,,\end{equation}
where $A_i\in \SL_2(\CC).$ It is well known, that the isomonodromic deformations of the system \eqref{fuchs} are described by the Schlesinger 
equations
\begin{equation}\label{schlesinger}
\frac{\partial A_i}{a_j}=\frac{[A_i,A_j]}{a_i-a_j}\quad 
{\rm if}\quad i\not=j, \quad {\rm and }\quad \frac{\partial A_i}{\partial
a_i}=-\sum_{j\not= i}\frac{[A_i,A_j]}{a_i-a_j}.\end{equation}

Set $A_4:=-(A_1+A_2+A_3)$ and let 
 $O_i$ denote the adjoint orbit, containing $A_i.$ 
Let $M$ be the quotient of 
$$\left\{ (A_1,\ldots,A_4)\in O_1\times \cdots \times O_4\mid
\sum A_i =0\right\}$$
modulo overall conjugation by $\SL_2(\CC),$ let 
$B:=\CC^3\setminus {\rm diagonals}$  and consider the 
trivial fibre bundle ${\cal M}:=M\times B \to B.$ 
Since the Schlesinger equations are invariant under 
overall conjugation,  they give rise to the {\em isomonodromy connection} $\nabla$ on 
the  fibre bundle ${\cal M}\to B.$ 

For $\CCC:=(C_1,\ldots,C_4),\,C_i:=\exp(2\pi \sqrt{-1}O_i)$ 
and $(a_1,a_2,a_3)\in B,$ consider 
the set 
$\Hom_\CCC(\pi_1(\CC\setminus  \{a_1,a_2,a_3\}),\SL_2(\CC))/\SL_2(\CC)$
where a simple (counterclockwise) loop around a point $a_i$ is mapped into $C_i$ and a simple (clockwise) loop around the set $\{a_1,a_2,a_3\}$
is mapped into $C_4.$ When $(a_1,a_2,a_3)$ varies in $B,$ these sets 
fit together to a fibre bundle $\psi:H\to B.$ This fibre bundle has a complete
flat connection $\Delta$ 
which is defined locally (say, $a_i \in D_i,\, i=1,2,3,$ 
where $D_i$ is a small open disk with center $a_i$),
by identifying representations which 
take the same values at simple loops $\gamma_i$ around the disks $D_i.$
For fixed $(a_1^o,a_2^o,a_3^o)\in B$ the set $\psi^{-1}(a_1^o,a_2^o,a_3^o)$ identifies
with the set 
$$ \NNN(\CCC):=\{(g_1,\ldots,g_4)\in C_1\times \cdots \times C_4\mid 
g_1\ldots g_4=1\}/\SL_2(\CC),$$
by fixing a homotopy base of $\pi_1(\CC\setminus \{a_1^o,a_2^o,a_3^o\}.$
 Let $[g_1,\ldots,g_4]\in \NNN(\CCC)$
 denote the equivalence class of $ (g_1,\ldots,g_4).$
By completeness, there is a unique global section of 
$\Delta$ through any point $P\in H.$ Moreover, the analytic 
continuation of a section through $[g_1,\ldots,g_4]\in \NNN(\CCC) $
along a closed path $\gamma$ at $ (a_1^o,a_2^o,a_3^o)$ is given by 
$[g_1,\ldots,g_4]^\gamma,$ where the (pure) braid group 
$B^3:=\pi_1(B,(a_1^o,a_2^o,a_3^o))$ acts as a subgroup 
of the full Artin braid group on 
tuples of matrices as in Thm. \ref{eigen}, (iv). 
Thus, sections of $\Delta$ with finite branching 
correspond to finite braid group orbits of the set $\NNN(\CCC).$\\

One has a natural map of fibre bundles 
$$ \nu:{\cal M} \longrightarrow H$$ 
which is induced by taking the monodromy of a Fuchsian system. By 
construction, the isomonodromy connection on ${\cal M}$ is the pullback
of $\nabla.$ Thus, finite branching 
sections of $\nabla$ correspond to finite braid group orbits on 
$\NNN(\CCC).$ \\

The inclusion 
$\iota: \CC\setminus \{0,1\}\hookrightarrow B,\,t \mapsto (0,t,1),$
induces an inclusion $\pi_1(\CC\setminus \{0,1\})\to B^3.$ Let $P^3$ be 
the image of this map.  
It is well known, that the Painlev\'e equation $\PVI$ arises from $\nabla$
by restricting the singularity positions to $(0,t,1),$ by a suitable choice 
of local coordinates of $M$ and by a parameter elimination. Thus, finite 
branching solutions of $\PVI$ are obtained from finite branching solutions 
of $\nabla,$ which correspond to finite braid group orbits on $\NNN(\CCC).$
It is expected, that these solutions of $\PVI$ are always algebraic, see 
\cite{Boalch03}, P. 1012.\\

\section{From triples pseudo-reflections in $\GL_3$ to triples 
in $\SL_2$}\label{boalch}

Let $G:=\SL_2(\CC).$ Any triple $(M_3,M_2,M_1) \in G^3$ gives rise
to a tuple ${\bf m}:={\bf m}(M_3,M_2,M_1):=(m_1,m_2,m_3,m_{13},m_{23},m_{13},m_{321})\in \CC^7,$
where 
$$ m_1:=\Tr(M_1),\quad m_2:=\Tr(M_2),\quad  m_3:=\Tr(M_3),$$
$$ m_{12}:=\Tr(M_1M_2),\quad  m_{23}:=\Tr(M_2M_3),\quad m_{13}:=\Tr(M_1M_3)$$
and $$m_{321}:=\Tr(M_3M_2M_1).$$
Such seven-tuples satisfy the so-called {\em Fricke relation}, see
\cite{Boalch04}, Section 2. Moreover, the induced 
braid group action can be explicitly determined, loc. cit..\\

Let $V:=\CC^3,$ $e_i\in V,\, \alpha_i\in V^*, \, i=1,2,3,$ and let 
$r_i:=1+e_i\otimes \alpha_i,\,
 i=1,2,3,$ be pseudo-reflections in $\GL_3(\CC).$ Choose complex numbers
$n_1,n_2,n_3,t_1,t_2,t_3$ such that $t_i$ is a choice of a square root 
of $\det(r_i),$ that the product $r_3r_2r_1$ has eigenvalues 
$\{n_1^2,n_2^2,n_3^2\}$ and that the square roots are chosen so that
$t_1t_2  t_3=n_1n_2n_3.$ Consider the following eight $\GL_3$-invariant 
functions on the space of triples of pseudo-reflections in $\GL_3:$
$$ t_1^2,\quad t_2^2,\quad t_3^2$$
$$ t_{12}:=\Tr(r_1r_2)-1,\quad t_{23}:=Tr(r_2r_3)-1,\quad 
t_{13}:=\Tr(r_1r_3)-1,$$
$$t_{321}:=n_1^2+n_2^2+n_3^2,\quad t_{321}':=(n_1n_2)^2+(n_2n_3)^2+(n_1n_3)^2.$$ Again, such eight-tuples satisfy certain relations and the braid group 
action can be explicitly determined (loc. cit.). \\

Now, it is shown in loc. cit. that the following map is invariant 
under the action of the braid group: Suppose that one is given 
a tuple 
$$ {\bf t}:={\bf t}(r_3,r_2,r_1):=(t_1,t_2,t_3,n_1,n_2,n_3,t_{12},t_{23},t_{13}),$$ 
associated to a triple of pseudo-reflections. Then one can define 
a map $\varphi$ taking ${\bf t}$ to some $\SL_2(\CC)$ data 
${\bf m}$ as follows: 
$$ m_1:=\frac{t_1}{n_1},\quad m_2:=\frac{t_2}{n_1}+\frac{n_1}{t_2},
\quad m_3:=\frac{t_3}{n_1}+\frac{n_1}{t_3},$$
\begin{equation} m_{12}:=\frac{t_{12}}{t_1t_2},\quad m_{23}:=\frac{t_{23}}{t_2t_3},
\quad m_{13}:=\frac{t_{13}}{t_1t_3},\end{equation}
$$ m_{321}:=\frac{n_2}{n_3}+\frac{n_3}{n_2}.$$

It is shown in loc.~cit., Thm1., that this map is $B_3$-invariant, thus 
finite $B_3$-orbits of $\SL_2(\CC)$-triples are obtained from 
triples of generators of three-dimensional complex reflection groups. These
in turn lead to new algebraic solutions, of $P_{VI}.$\\

The following result shows, how the map $\varphi $ is related 
to the middle convolution:

\begin{thm} \label{phi} Let ${\bf r}:=(r_3,r_2,r_1)\in 
\GL_3(\CC)$ be a triple of pseudo-reflections and let 
$t_1,\ldots,n_1,\ldots$ and ${\bf t}$ be as above. Suppose that 
$r_3r_2r_1$ has three distinct eigenvalues. 
Let $\lambda\in \CC^\times$ be an eigenvalue of $r_3r_2r_1,$ 
$$(M_3',M_2',M_1'):=MC_\lambda({\bf r})$$ and 
$$ (M_3,M_2,M_1):= (n_3M_3',n_2M_2',n_1M_1').$$
Then 
$$ {\bf m}(M_3,M_2,M_1)=\varphi({\bf t}).$$
\end{thm}

\proof

\section{Some examples}\label{examples}

\Endproof

\bibliographystyle{plain}
\bibliography{p}

Michael Dettweiler

IWR, Universit\"at Heidelberg,

INF 368

69121 Heidelberg, Deutschland

e-mail: michael.dettweiler@iwr.uni-heidelberg.de\\

Stefan Reiter

Technische Universit\"at Darmstadt

Fachbereich Mathematik AG 2

Schlo\ss gartenstr. 7

64289 Darmstadt, Deutschland

email: reiter@mathematik.tu-darmstadt.de

\end{document}